\input amstex
\documentstyle{amsppt}
\pageheight{50.5pc}
\pagewidth{32pc}
\define\1{\hbox{\rm 1}\!\! @,@,@,@,@,\hbox{\rm I}}

\topmatter
\title
{
Limit theorems for number of diffusion processes which did not absorb
by boundaries.
\footnotemark{}
}
\endtitle
\footnotetext{This research was supported (in part) by the
Ministry of Education and Science of Ukraine, project
No 01.07/103 and University Salerno, Italy.}
\author
{Aniello Fedullo and Vitalii A. Gasanenko}
\endauthor

\address
{Universita Degli Studi Di Salerno, Via Ponte Don Melillo, 84084
 Fisciano (SA) Italia}
\endaddress
\email
{afedullo\@unisa.it}
\endemail
\address
{Institute of Mathematics,National Academy of
Science of Ukraine, Tereshchenkivska 3, 01601, Kiev, Ukraine}
\endaddress
\email
{gs\@imath.kiev.ua}
\endemail
\subjclass
{60 J 60}
\endsubjclass
\endkeywords
\subjclass
{60 J 60}
\endsubjclass
\email
{afedulo\@unisa.it}
\endemail
\email
{gs\@imath.kiev.ua}
\endemail
\abstract
We have random number of independent diffusion  processes
 with absorption on boundaries in some region  at initial time $t=0$.
 The initial numbers and positions of processes in region
is defined by Poisson random measure.
It is required to estimate of number
of the  unabsorbed processes for the fixed time ~$\tau>0$.
The Poisson random measure depends on $\tau$ and $\tau\to\infty$.
\endabstract

\endtopmatter

\document

Consider the  set of  independent random  diffusion processes ~~~$\xi_{k}(t)
,\quad k=\overline{1,N},$

$t\geq 0,~~\xi_{k}(0)=x_{k},~~
x_{k}\in Q\subset R^{d}$.
We wish to investigate of distribution
 of the number of the  processes ~$\xi_{k}(t)$~
which was into $Q$ for all moments of time  $t\leq \tau$.

Let domain  ~~$ Q\subset R^{d}$~~ be  open connected region and it is limited by
smooth surface $\partial Q$. All processes ~$\xi_{k}(t)$~
are diffusion processes with absorption on the boundary ~$\partial Q$.
 These processes are solutions
of the following stochastic differential  equations in $Q$

$$
d\xi(t)=a(t,\xi(t))dt + \sum\limits_{i=1}^{d}b_{i}(t,\xi (t))dw^{(k)}_{i}(t)
\quad \xi(t)\in R^{d}\eqno(1)
$$

$$
 b_{i}(t,x),~ a(t,x): R_{+}\times R^{d}\to R^{d}.
$$

with an initial condition: ~$\xi(0)=x_{k}\in D.$

Here the ~$W^{(k)}(t)=(w_{i}^{(k)}(t),\quad 1\leq i\leq d),\quad 1\leq k\leq N$~
 are independent in totality $d$- dimensional Wiener processes.

Thus, these processes have the identical diffusion matrices and
shift vectors , but they have different initial states.

Let $Q$ is bounded and boundary $\partial Q$ is  Lyapunov surface
~$C^{(1,\lambda)}$.

 The initial number and positions of processes
 are defined  by the random Poisson measure $\mu(\cdot,\tau)$ in $Q$:
$$
P(\mu(A,\tau)=k)=\frac{m^{k}(A,\tau)}{k!}e^{-m(A,\tau)}
$$
where~$m(\cdot,\tau)$ is finitely  additive positive mesure on ~$Q$
for fixed $\tau$.

This task was offered in [1] as the mathematical model of
 practice problem.
The authores in article [2] investigated case when
 initial number and positions of diffusion processes
 are defined  by determinate limited measure
~$N(B,\tau)$.
Where the ~$N(B,\tau)$~ is equal to number of points ~$x_{k}$~ in a set ~$B$~
and ~$N=N(Q,\tau)<\infty$ for fixed $\tau>0$.

We consider the following case

$a(t,x)=a=(\underbrace{0,\dots,0}_{d}),\quad b_{i}(t,x)=b_{i}=
(b_{i1},\dots,b_{id}),~~1\leq i\leq d;$

We  define  matrix
 ~$\sigma=B^{T}B,\quad B=(b_{ij}),~1\leq i,j\leq d$

$\sigma=(\sigma_{ij}),1\leq i,j\leq d$ and differential operator
$A:\sum\limits_{1\leq i,j\leq d}\sigma_{ij}\frac{\partial^{2}}
{\partial x_{i}\partial x_{j}}.$

Let $\sigma$ be a matrix with  the following property

$$
\sum\limits_{1\leq i,j\leq d}\sigma_{ij}z_{i}z_{j}
\geq \mu |\vec z|^{2}.
$$

Here ~$\mu$,~ there is  fixed positive number, and
$\vec z=(z_{1},\cdots, z_{d})$~ there is an arbitrary real vector.

This operator acts in the following space

$$
H_{A}=\{u: u\in L_{2}(Q)\cap Au\in L_{2}(Q)\cap u(\partial Q)=0\}
$$

with inner product $(u,v)_{A}=(Au,v)$.Here $(,)$ is inner product in $L_{2}(Q)$.
 The operator ~$A$~ is  positive operator.

 It is known [3] that the following eigenvalues problem

$$
Au=-\lambda u,\quad u(\partial Q)=0
$$

has infinity set of real eigenvalues $ \lambda_{i}\to\infty$ and

$$
0<\lambda_{1}<\lambda_{2}<\cdots<\lambda_{s}<\cdots.
$$

The corresponding eigenfunctions

$$
f_{11},\dots,f_{1n_{1}},\cdots,f_{s1},\dots,f_{sn_{s}},\cdots
$$

form complete system of functions both in $H_{A}$ and

$L_{2}^{0}(Q):=
\{u: u\in L_{2}(Q)\cap u(\partial Q)=0\}$. Here the number $n_{k}$ is equal to
multiplicity  of eigenvalue  $\lambda_{k}$.

We  denote by ~$\eta(\tau)$~ the number of remaining processes in the region ~$D$~
at time instant ~$\tau$.

We also assume that $\sigma$-additive  measure $\nu$ is given on the
$\Sigma_{\nu}$- algebra sets of 

$Q,\quad \nu(Q)<\infty.$
All eigenfunctions
$f_{ij}:Q\to R^{1}$ and all measures $m(\cdot,\tau)$~ are $(\Sigma_{\nu},\Sigma_{Y})$ measurable. Here
$\Sigma_{Y}$ is system of Borel sets of $R^{1}$.
Let $\Rightarrow$ denotes the  weak convergence of
random values.

Put

$$
g(\tau)=\exp\left(-\frac{\tau}{2}\lambda_{1}\right).
$$

\proclaim
{\bf Theorem 1}{\it We suppose that $m(\cdot, \tau)$ satisfies the condition

$$
\lim\limits_{\tau\to\infty}
m(B,\tau)g(\tau)= \nu(B),\quad B\in \Sigma_{\nu}.
$$

Then $\eta(\tau)\Rightarrow\eta$ if $\tau\to\infty$ where $\eta$ has Poisson
distribution function with the

parameter $a=\int\limits_{Q} F(x)d\nu$ and
$F(x)=\sum\limits_{i=1}^{n_{1}}f_{1i}(x)c_{1i},\quad
c_{1i}=\int\limits_{Q}f_{1i}(\bar x)d\bar x.$}
\endproclaim
\demo
{\bf Proof} We consider the following initial-boundary problem

$$
\frac{\partial u}{\partial t}=\frac{1}{2}\sum\limits_{1\leq i,j\leq d}
\sigma_{ij}\frac{\partial^{2}u}{\partial x_{i}\partial x_{j}}\quad x\in Q;
$$

$$
u(0,x)= 1\quad\hbox{if}\quad x\in Q;
$$

$$
u(t,x)=0 \quad \hbox{if}\quad x\in\partial Q,~~~t\geq 0\eqno(2)
$$

It is known [4], that ~ $u(\tau, x)$~ is equal to probability of remaining
 in the
 region ~$Q$~  at time instant ~$\tau$~ of a diffusion process from (1) which
 occurs at the point $(0,x)$~at the  initial moment
(~$\xi(0)= x,~~  x\in Q$).
We  designate through $\gamma_{k}=(x^{k}_{1}, \cdots ,x_{d}^{k})$
the initial position of
$k$-th process. We define the value of $u(\tau,\gamma_{k})$.

We  define a particular solution of (2) in form

$$
u(t, x)= u_{1}(t)u_{2}(x).
$$

The ordinary argumentaion leads to  definition of
joined constant $\lambda$:

$$
2\frac{1}{u_{1}}\frac{\partial u_{1}}{\partial t}
=\frac{A u_{2}}{u_{2}}=- \lambda.
$$

We obtain the following system of tasks due the latter one

$$
A u_{2}= -\lambda u_{2};\quad u_{2}(\partial Q)=0.\eqno(3)
$$

$$
\frac{\partial u_{1}}{\partial t}=
-\frac{\lambda}{2} u_{1};
\quad
u_{1}(0)= 1\eqno(4),
$$

It is clear that $u_{1}(t,\lambda)=
\exp(-\frac{t}{2}\lambda )$
is solution
of (5) .
The soluton of (3) was described above. We  assume
that system of functions $\{f_{ij}( x), i\geq 1, 1\leq j\leq n_{i}\}$
is orthonormalized with respect to space $L_{2}^{0}(Q)$.

The general solution of problem (2) has the following form

$$
u(t,x)=\sum\limits_{j=1}^{\infty}
\exp(-\frac{t}{2}\lambda_{j} )
\sum\limits_{m=1}^{n_{j}}c_{jm}f_{jm}( x)
$$

where coefficients $c_{jm}$ are equal to coefficients of decomposition
of initial value (unit)  by system of functions $f_{jm}$:
$c_{jm}=\int\limits_{Q}f_{jm}( x)d x$. The Parseval - Steklov equality
is true for these coefficients:

$$
\sum\limits_{j=1}^{\infty}\sum\limits_{m=1}^{n_{j}}c_{jm}^{2}=|Q|.\eqno(5)
$$

Put $F( x)=\sum\limits_{m=1}^{n_{1}}c_{1i}f_{1m}( x)$. The
function ~$F(x)$~ is continuous and bounded function on the $\bar Q$. Since
$u(t,x)$ is probability,it is not diffucult to show that $F( x)\geq 0$ for all $x\in Q$.
Let $M=\sup\limits_{ x\in Q}F(x)$.
We introduce the following sets

$$
B_{k,n}=\{  x\in Q:  \frac{Mk}{n}< F( x)\leq\frac{M(k+1)}{n} \}
$$

Here ~$0\leq k\leq n-1$~  and $n>1$.

Let us denote by
$\zeta_{k,n}(\tau), ~1\leq~k\leq n$~ the number of unabsorbed processes
at time instant ~$\tau$ which occur in the region ~$B_{k,n}$ at the initial time.
These values are independent in totality by assumption.
The distribution function of  $\zeta_{k,n}(\tau)$~
is defined by the following formula

$$
P(\zeta_{k,n}(\tau)=l)
=\sum\limits_{d=l}^{\infty}P(\mu(B_{k,n},\tau)=d)\times
$$
$$
\times\sum\limits_{1\leq i_{1},\cdots,i_{l}\leq d,i_{m}\ne i_{j},m\ne j}
\prod\limits_{k=1}^{l}u(\tau,\gamma_{i_{k}})
\prod\limits^{d}_{s=l+1, i_{s}\notin (i_{1},\cdots,i_{l}),i_{m}\ne i_{j}}
(1- u(\tau,\gamma_{i_{s}})),\quad l=0,1,\dots.
$$

Here ~$x_{i_{j}}\in B_{k,n}$.

We set

$$
a_{k,n}(\tau)=\min\limits_{ x\in B_{k,n}}u(\tau, x),\quad
\bar a_{k,n}(\tau)=1- a_{k,n}(\tau);
$$
$$
b_{k,n}(\tau)=\max\limits_{ x\in B_{k,n}}u(\tau, x),\quad
\bar b_{k,n}(\tau)=1- b_{k,n}(\tau).
$$

Now

$$
J_{k,n}(l,\tau):=
\sum\limits_{d=l}^{\infty}\frac{m^{d}(B_{k,n},\tau)}{d!}
\exp(- m(B_{k,n},\tau))
 C_{d}^{l}a^{l}_{k,n}(\tau)\bar b^{d-l}_{k,n}(\tau)\leq
$$
$$
\leq P(\zeta_{k,n}(\tau)=l)\leq
$$
$$
\sum\limits_{d=l}^{\infty}
\frac{m^{d}(B_{k,n},\tau)}{d!}\exp(- m(B_{k,n},\tau))
 C_{d}^{l}b^{l}_{k,n}(\tau)\bar a^{d-l}_{k,n}(\tau)=:
I_{k,n}(l,\tau).\eqno(6)
$$

Further

$$
J_{k,n}(l,\tau)=
\frac{\left(m(B_{k,n},\tau)a_{k,n}(\tau)\right)^{l}}{l!}
\exp(- m(B_{k,n},\tau))
\sum\limits_{d=l}^{\infty}
\frac{\left(\bar b_{k,n}(\tau)m(B_{k,n},\tau)\right)^{d-l}}{(d-l)!}
=
$$
$$
=\frac{\left(m(B_{k,n},\tau)a_{k,n}(\tau)\right)^{l}}{l!}
\exp(- b_{k,n}(\tau)m(B_{k,n},\tau));
$$

By analogy:

$$
I_{k,n}(l,\tau)=
\frac{\left(m(B_{k,n},\tau)b_{k,n}(\tau)\right)^{l}}{l!}
\exp(- a_{k,n}(\tau)m(B_{k,n},\tau)).\eqno(7)
$$

We introduce the following generating functions

$$
\varphi(\tau,s)=\sum\limits_{l\geq 0}s^{l}P(\eta(\tau)=l).
$$
$$
\varphi_{k,n}(\tau,s)=\sum\limits_{l\geq 0}s^{l}P(\zeta_{k,n}(\tau)=l),
\quad k=\overline{0,n-1},\quad 0\leq s\leq 1.
$$

By the construction, $\eta(\tau)$~ can be represented in the form
 ~$\eta(\tau)=\zeta_{1,n}+\cdots+\zeta_{n-1,n}(\tau)$.

Thus

$$
\varphi(\tau,s)=\prod_{k=0}^{n-1}\varphi_{k,n}(\tau,s).\eqno(8)
$$

Combining  (6)-(8), we conclude that

$$
\exp\{(sa_{k,n}(\tau)-b_{k,n}(\tau))m(B_{k,n},\tau)\}\leq
\varphi_{k.n}(\tau,s)\leq
$$
$$
\leq \exp\{(sb_{k,n}(\tau)-a_{k,n}(\tau))m(B_{k,n},\tau)\}
$$

and

$$
\exp\left\{\sum\limits_{k=0}^{n-1}(sa_{k,n}(\tau)-b_{k,n}(\tau))
m(B_{k,n},\tau)\right\}\leq
\varphi(\tau,s)\leq
$$
$$
\leq
\exp\{\sum\limits_{k=0}^{n-1}(sb_{k,n}(\tau)-a_{k,n}(\tau))
m(B_{k,n},\tau)\}.\eqno(9)
$$

Since function ~$u(\tau,x)$~ is continuous function in
$ x\in Q$, there exit a points $x_{*},~~x^{*}\in \bar B_{k,n}$
such that
the following equalities have place

$$
a_{k,n}(\tau)=
u_{1}(\tau,\lambda_{1})F(x_{*})+
\sum\limits_{k\geq 2}u_{1}(\tau,\lambda_{k})\sum\limits_{m=1}^{n_{k}}c_{km}
f_{km}(x_{*}),
$$
$$
b_{k,n}(\tau)=
u_{1}(\tau,\lambda_{1})F( x^{*})+
\sum\limits_{k\geq 2}u_{1}(\tau,\lambda_{k})\sum\limits_{m=1}^{n_{k}}c_{km}
f_{km}(x^{*}),
$$

here ~$ x_{*}:=x_{*}(k,n,\tau),~~x^{*}:=x^{*}(k,n,\tau)$.

Now, we can rewrite  the sums in  exponentes from (9) in the following forms

$$
\sum\limits_{k=0}^{n-1}
\left(sF( x_{*})- F(x^{*})\right)
\exp(-\frac{\tau}{2}\lambda_{1} )
m(B_{k,n},\tau)+
$$
$$
+
\sum\limits_{k=0}^{n-1}
\exp(-\frac{\tau}{2}\lambda_{1} )
m(B_{k,n},\tau)
\sum\limits_{j\geq 2}
\exp\left(-\frac{\tau}{2}(\lambda_{j}-\lambda_{1})\right)
\sum\limits_{m=1}^{n_{j}}c_{jm}
(sf_{jm}( x_{*})-f_{jm}( x^{*})) , \eqno(10)
$$

$$
\sum\limits_{k=0}^{n-1}(sF( x^{*})- F( x_{*}))
\exp(-\frac{\tau}{2}\lambda_{1} )
m(B_{k,n},\tau)+
$$
$$
+\sum\limits_{k=0}^{n-1}
\exp(-\frac{\tau}{2}\lambda_{1} )
m(B_{k,n},\tau)
\sum\limits_{j\geq 2}
\exp\left(-\frac{\tau}{2}(\lambda_{j}-\lambda_{1})\right)
\sum\limits_{m=1}^{n_{j}}c_{jm}
(sf_{jm}( x^{*})- f_{jm}( x_{*})), \eqno(11)
$$

We calculate limit of (10) if ~$\tau\to\infty$.
The first sum of (10) convergers
to the following limit
under the condition of theorem

$$
\sum\limits_{k=0}^{n-1}sF( x_{*})\nu(B_{k.n})
-\sum\limits_{k=0}^{n-1}F( x^{*})\nu(B_{k.n}).
$$

This is difference of two integral sums which
has the following limit under ~$n\to\infty$ (see [5])
$$
(s-1)\int\limits_{Q}F(x)\nu(d x).
$$

Put
$$
s_{\tau}( x)=
\sum\limits_{j\geq 2}\exp\left(-\frac{\tau}{2}(\lambda_{j}-\lambda_{1})\right)
\sum\limits_{m=1}^{n_{j}}c_{km}
f_{km}(x).
$$

We consider  sums of  eigenfunctions in the form

$$
e(x,\lambda)=\sum\limits_{\lambda_{j}\leq \lambda}f^{2}_{jl}(x)
$$

The following result is proved in monography [6,Thm. 17.5.3]
$$
\sup\limits_{x\in Q}e(x,\lambda)\leq C \lambda^{\frac{d}{2}}.
$$

Asymptotic characteristic of eigenvalues ~$\lambda_{j}$~ under
~$j\to\infty$~ is defined by the following inequalities [3, sec.
18]
$$
c_{1}j^\frac{2}{d}\leq \lambda_{j}\leq c_{2} j^\frac{2}{d}, \quad
\hbox{where}\quad c_{1},~c_{2}=const.
$$

The latter one, (5) and  Caushy-Bunyakovskii inequality
lead to the following convergence under
$\tau\to\infty$

$$
|s_{\tau}( x)|\leq
\sum\limits_{j\geq 2}\exp\left(-\frac{\tau}{2}(\lambda_{j}-\lambda_{1})\right)
\sqrt{\sum\limits_{m=1}^{n_{j}}
c^{2}_{jm}}\sqrt{\sum\limits_{m=1}^{n_{j}}f_{jm}^{2}(x)}\leq
$$
$$
\leq \sqrt{C}\sum\limits_{j\geq 2}\lambda_{j}^{\frac{d}{2}}\exp(-\frac{\tau}{2}
(\lambda_{j}-\lambda))\sqrt{\sum\limits_{m=1}^{n_{j}}c_{jm}^{2}}\leq
$$
$$
\leq \sqrt{C} \sqrt{\sum\limits_{j\geq 2} \lambda_{j}^{d}
\exp(-\tau (\lambda_{j}-\lambda_{1}))}
\sqrt{\sum\limits_{j\geq 2}\sum\limits_{m=1}^{n_{j}}c^{2}_{jm}}\to 0.
$$

Thus the second sum from(10) convergences to zero.

The similar considerations apply to (11). Proof is complete.
\enddemo

{\it\bf  Example.}
Now we  apply the general approach to the particular case.

We consider the case if  $Q$ is circle
~~$Q=\{(x,y): x^{2}+y^{2}\leq r_{0}^{2}\}$.
We assume that the diffusion processes occurs at the  point
$(x_{k},y_{k})\in Q$ at the initial time.

The processes are described in ~$Q$~
by the following stochastic differential equations

$$
d\xi(t)=
\sum\limits_{i=1}^{2}b_{i}dw_{i}(t)\eqno(12)
$$
$$
\xi(0)=\xi_{0}=(x_{k},y_{k}).
$$

where $ b_{1}=(\sigma,0),b_{2}=(0,\sigma) $ and

$W(t)=(w_{i}(t),i=1,2)$ is  2-dimensional Wiener process.

We assume that the equation (12) defines a diffusion process with absorption
on the boundary
 $\partial Q=\{(x,y,z): x^{2}+y^{2}=r_{0}^{2}\}$.

In follows that the $J_{0}(x), J_{1}(x)$ are Bessel functions zero and first order.
They are defined as the solutions of the following  equations

$$
\frac{d^{2} y}{d x^{2}} + \frac{1}{x}\frac{d y}{d x} + (1-\frac{n^{2}}{x^{2}})
=0,
$$
$$
y(x_{0})=0,~~~  (x_{0}=\sqrt{\lambda}r);\quad  |y(0)|<\infty;
$$

under $n=0$ and $n=1$.

The value of $\mu_{m}^{(0)}$ is equal to $m$- th
root of the equation $J_{0}(\mu)=0$ [7,8].

Let $mes(\cdot)$  denotes the Lebesgue measure.

We set

$$
f(\tau):=\exp\left(-\frac{\tau}{2}
\Bigl(\frac{\sigma \mu_{1}^{(0)}}{r}\Bigr)^{2}
\right).
$$

\bigskip

We suppose that $m(\cdot,\tau)$ holds the condition
$$
m(\cdot,\tau)f(\tau) \Rightarrow mes(\cdot) \quad \hbox{if}\quad \tau\to\infty.
$$

In this case the system of tasks (3),(4) has the following form

$$
\triangle u_{2}= -\mu u_{2}, ~~(x,y)\in C ;\quad
u_{2}(x,y)=0 \quad \hbox{if} \quad x^{2}+y^{2}=r_{0}^{2}\eqno(13),
$$
$$
\frac{\partial u_{1}}{\partial t}=-\frac{\sigma^{2}}{2} \mu u_{1}, \quad
u_{1}(0)=1.\eqno(14)
$$

According to general approach for construction of solution ~$u(t,x,y)$~
(see,for example, [7, sec.1V ]
we rewrite the task of (13) in polar coordinates:
 $u_{3}(r,\varphi):=
u_{2}(r\cos\varphi,r\sin\varphi)$.
The $u_{3}$ is solution the following
problem

$$
\frac{\partial^{2} u_{3}}{\partial r^{2}}+\frac{1}{r}\frac{\partial u_{3}}
{\partial r}+ \frac{1}{r^{2}}\frac{\partial^{2}u_{3}}{\partial \varphi^{2}}
+\mu u_{3}=0,
$$

$$
u_{3}(r_{0},\varphi)=0.
$$

We obtain

$$
u(t,x,y)=u(t,r)=\sum\limits_{m=1}^{\infty}c_{m}
J_{0}\left(\frac{\mu_{m}^{(0)}}{r_{0}}r\right)
\exp\left(-\frac{t}{2}
\left(\frac{\sigma\mu^{0}_{m}}{r_{0}}\right)^{2}\right),
$$

where $c_{m}=2\left(m_{m}^{(0)}J_{1}(\mu_{m}^{(0)})\right)^{-1}$.

The function
~$J_{0}\left(\frac{\mu_{1}^{(0)}}{r_{0}}r\right)$~ is strictly decreasing
function if ~$0\leq r\leq r_{0}$~. Thus we can construct the
partitions ~$B_{k,n}$~ by the following partitions
$$
\tilde B_{k,n}=\left\{ (x,y)\in C: \frac{r_{0}k}{n}< \sqrt{x^{2}+y^{2}}
\leq \frac{r_{0}(k+1)}{n}\right\},\quad 0\leq i\leq n-1.
$$

Now ~$mes(\tilde B_{k,n})=g(\frac{k+1}{n})-g(\frac{k}{n})$~ where ~
~$g(x)=\pi r_{0}^{2}x^{2}, ~~~0\leq x\leq 1$.

Finally, the parameter of Poisson distribution is equal to

$$
a=2\left(m_{1}^{(0)}J_{1}(\mu_{1}^{(0)})\right)^{-1}2\pi r_{0}^{2}
\int\limits_{0}^{1}J_{0}(\mu_{1}^{(0)}x)xdx=\pi
\left(\frac{2 r}{\mu_{1}^{(0)}}\right)^{2}.
$$

We used the following known relation
~~$\alpha J_{0}(\alpha)=\left[\alpha J_{1}(\alpha)\right]'$ [7,p.466]
for calculation of the latter integral.

\enddocument